\documentclass[10pt,letterpaper]{article}
\usepackage{fullpage}

\usepackage{defs}
\usepackage{words}
\usepackage{art-thm}

\title{A hyperelliptic Hodge integral}
\author{Jonathan Wise}
\date{\today}

\begin{document}

\maketitle

\section{Introduction}

We work over $\CC$.  We will use the theories of orbifold stable maps and orbifold Gromov--Witten theory as developed in~\cite{AV} and~\cite{AGV}.  Our notation for the moduli space of degree $\beta$ orbifold stable maps with $n_1$ ordinary marked points and $n_2$ orbifold points will be
\begin{equation*}
\overline{M}(X; n_1, n_2; \beta) ,
\end{equation*}
which is an open substack of the corresponding Artin stack of pre-stable maps, ${\frk M}(X; n_1, n_2; \beta)$.  (Since we will only deal with $\ZZ / 2 \ZZ$ stabilizers in this note, we will not need a more detailed notation.)

  Let $C$ be the universal curve over ${\frk M} = {\frk M}(B(\ZZ / 2 \ZZ); 0, 2g + 2)$.  For each $i = 1, \ldots, 2g + 2$, there is a closed substack $D_i$, the $i$-th universal $\ZZ / 2 \ZZ$-gerbe over ${\frk M}$.  Let $N_{D_i / C}$ be the normal bundle of $D_i$ in $C$ and define $L_i$ to be the line bundle 
\begin{equation*}
L_i = \pi_\ast \left( N_{D_i / C}^\vee \tensor \rho_1 \right) .
\end{equation*}
where $\rho_1$ is the non-trivial representation of $\ZZ / 2 \ZZ$, viewed as a line bundle on $D$ pulled back from $B(\ZZ / 2 \ZZ)$.

\begin{remark}
Our definition of $L_i$ coincides with the cotangent line bundle on the universal hyperelliptic curve over ${\frk M}(B(\ZZ / 2 \ZZ); 0, 2g + 2)$.  Indeed, let $p : \tilde{C} \rightarrow C$ be the base change of $(\rm point) \rightarrow B(\ZZ / 2 \ZZ)$ via the universal map $C \rightarrow B(\ZZ / 2 \ZZ)$ and define $\tilde{D}_i$ analogously.  Then $\pi_\ast (N_{D_i / C} \tensor \rho_1)$ can be identified with the pushforward via $\pi p$ of the $-1$-eigenspace of $p^\ast N_{D_i / C}^\vee = N_{\tilde{D}_i / \tilde{C}_i}^\vee$.  Since the hyperelliptic involution acts nontrivially on the fiber of the cotangent bundle at a Weierstrass point, this is just $N_{\tilde{D}_i / \tilde{C}_i}$ which is the usual cotangent line bundle.
\end{remark}

In view of the remark, it is legitimate to say $c_1(L_i) = \psi_i$.

We also have the hyperelliptic Hodge bundle, whose dual is defined to be
\begin{equation*}
{\bf E}^\vee = R^1 \pi_\ast (\rho_1)
\end{equation*}
where $\pi$ is the map from the universal curve $C$ to $\overline{M}(B(\ZZ / 2 \ZZ), 2g + 2)$.

\begin{remark}
This definition of the Hodge bundle coincides with the usual definition as $R^1 (\pi p)_\ast {\cal O}_{\tilde{C}}$ (where $p : \tilde{C} \rightarrow C$ is defined as in the last remark).  Indeed, we can identify $R^1 (\pi p)_\ast {\cal O}_{\tilde{C}}$ since $p_\ast {\cal O}_{\tilde{C}} \cong {\cal O}_C \oplus \left( {\cal O}_C \tensor \rho_1 \right)$ and ${\cal O}_C$ has no higher cohomology (because $C$ has genus $1$).
\end{remark}

It is therefore justified to write $c_i({\bf E}) = \lambda_i$.

\begin{theorem-in-sect}
We have
\begin{equation*}
\int_{\overline{M}(B(\ZZ / 2 \ZZ), 2g + 2)} \frac{c({\bf E}^\vee)^2}{c(L_1^\vee)} = \int_{\overline{\cal H}_g} \frac{ (1 - \lambda_1 + \cdots + (-1)^g \lambda_g )^2 }{1 - \psi_1} = \left( - \frac{1}{4} \right)^g .
\end{equation*}
\end{theorem-in-sect}

The first equality was proved in the two remarks above.  The second equality will be proven by interpreting the integral as a Gromov--Witten invariant on the weighted projective space (Section~\ref{sect:interp}) and evaluating it recursively using the WDVV equations (Section~\ref{sect:eval}).

The application for this calculation is~\cite{JW}, where it is used to relate the genus zero Gromov--Witten invariants of $[\Sym^2 \PP^2]$ and the enumerative geometry of hyperelliptic curves in $\PP^2$.

\section{A Gromov--Witten invariant of $\PP(1,1,2)$}
\label{sect:eval}

Let $\overline{M}(\PP(1,1,2); n_1, n_2; \beta)$ be the moduli space of genus zero orbifold stable maps to $\PP(1,1,2)$ with $n_1$ ordinary marked points and $n_2$ orbifold marked points and degree $\beta$.  The degree is evaluated by integrating $c_1({\cal O}(1))$ over the curve and so is an element of $\frac{1}{2} \ZZ$.

The virtual dimension is given by the formula
\begin{equation*}
\vdim \overline{M}(\PP(1,1,2); n_1, n_2; \beta) = \dim \PP(1,1,2) - 3 + \int_\beta c_1(T \PP(1,1,2)) + n_1 + n_2 - \sum_{i = 1}^{n_2} \age(x_i)
\end{equation*}
where $x_i$, $i = 1, \ldots, n_2$ is the set of orbifold marked points and $\age(x_i)$ is the sum of the $t_j$ such that the eigenvalues of the action of the stabilizer of $x_i$ acting on $T \PP(1,1,2)$ are $e^{2 \pi i t_j}$, $j =1, \ldots, n_2$, listed with multiplicity.  If $f : C \rightarrow \PP(1,1,2)$ is a representable map then any orbifold point of $C$ must be carried by $f$ to the unique stacky point of $C$, which is represented by $(0,0,1)$.  The automorphisms act with eigenvalues $-1,-1$ on the fiber of the tangent bundle at this point, so the age is $1$.

The Euler sequence here is
\begin{equation*}
0 \rightarrow {\cal O} \rightarrow {\cal O}(1) \oplus {\cal O}(1) \oplus {\cal O}(2) \rightarrow T \PP(1,1,2) \rightarrow 0
\end{equation*}
so $c_1(T \PP(1,1,2)) = 4 c_1({\cal O}(1))$.  Thus,
\begin{equation*}
\vdim \overline{M}(\PP(1,1,2); n_1, n_2; \beta) = 4d - 1 + n_1
\end{equation*}
where $d = \int_\beta c_1({\cal O}(1))$.

The inertia stack of $\PP(1,1,2)$ is $\PP(1,1,2) \amalg B (\ZZ / 2 \ZZ)$ and the rigidified inertia stack is $\PP(1,1,2) \amalg ({\rm point})$.  Let $\gamma$ be the fundamental class of the second component.  Let $p$ be the class of an ordinary point in $\PP(1,1,2)$.  We'll compute the invariant,
\begin{equation*}
\GW{p, \gamma, \ldots, \gamma}{\frac{1}{2}} .
\end{equation*}
Let's put $h = c_1({\cal O}(1))$.  Then $p = 2 h^2$.  

A schematic of the orbifold Chow ring  of $\PP(1,1,2)$ with its structure as a graded vector space is shown below.
\begin{equation*} \begin{array}{cc|c|c}
 & 0 & 1 & 2 \\
\PP(1,1,2) & \QQ & \QQ h & \QQ p \\
({\rm point}) & & \QQ \gamma
\end{array} \end{equation*}
It is easy to see that $\overline{M}(\PP(1,1,2); 1, 2; 0) \cong B(\ZZ / 2 \ZZ)$ and therefore that $\gamma^2 = \frac{1}{2} p = h^2$.  Therefore a presentation of the orbifold Chow ring is $\QQ[h, \gamma] / (h^2 - \gamma^2)$.  Note in particular that this satisfies Poincar\'e duality.

\begin{lem-in-sect} \label{gw-props}
The Gromov--Witten invariants of $\PP(1,1,2)$ have the following properties.
\renewcommand{\labelenumi}{(\alph{enumi})}
\begin{enumerate}
\item \label{gw-prop-1} If $\GW{\gamma^{\tensor n}, \alpha}{0} \not= 0$ then $n = 2$ and $\alpha = 1$.
\item \label{gw-prop-2} The invariant $\GW{\gamma^{\tensor n}, h, \ast}{0}$ is zero for all $n$.
\end{enumerate}
\end{lem-in-sect}
\begin{proof}
For (a), the invariant is computed on the moduli space $\overline{M}(\PP(1,1,2); a, n; 0)$ which has virtual dimension $a - 1$.  Hence the invariant will be zero unless $a = 1$ (so $\alpha$ comes from the untwisted sector) and $\alpha = 1$.  But then the invariant will be zero by the unit axiom unless $n = 2$.

For (b), it is sufficient by linearity to show that $\GW{\gamma^{\tensor n}, 2h, \ast}{0} = 0$.  But the Chow class $2h$ can be represented by a line that doesn't pass through the uniqe orbifold point $(0,0,1)$.  Since this is a degree zero invariant, this means it is computed on an empty moduli space, i.e., it is zero.
\end{proof}

The WDVV equations give
\begin{equation*}
\sum_{\substack{a + b = 2g - 1 \\ d_1 + d_2 = \frac{1}{2}}} \GW{ \GW{h, h, \gamma^{\tensor a}, \ast}{d_1}, \gamma, \gamma, \gamma^{\tensor b} }{d_2} = \sum_{\substack{a + b = 2g - 1\\ d_1 + d_2 = \frac{1}{2}}} \GW{ \GW{ h, \gamma, \gamma^{\tensor a}, \ast}{d_1}, h, \gamma, \gamma^{\tensor b} }{d_2} .
\end{equation*}
where $d_1$ and $d_2$ can take the values $0$ and $\frac{1}{2}$ in the sumes.

Consider first the right side of the equality.  One of the $d_i$ must be zero, so consider the invariant $\GW{h, \gamma, \ldots, \gamma, \ast}{0}$.  This is zero by Lemma~\ref{gw-props}~(b).  On the left side, note that if $d_1 = 0$ then the corresponding term of the sum will be zero by the divisor axiom unless $a = 0$ also.  Thus we get
\begin{equation*}
\GW{h^2, \gamma^{\tensor (2g + 1)}}{\frac{1}{2}} + \frac{1}{4} \sum_{a + b = 2g - 1} \GW{ \GW{ \gamma^{\tensor a}, \ast }{\frac{1}{4}}, \gamma^{\tensor (b + 2)}}{0} = 0 .
\end{equation*}
But $\GW{\gamma^{\tensor n}, \ast}{0} = 0$ for $n > 2$ by Lemma~\ref{gw-props}~(a).  Thus we are left with
\begin{equation*}
\GW{h^2, \gamma^{\tensor (2g + 1)}}{\frac{1}{2}} = - \frac{1}{4} \GW{\gamma^2, \gamma^{\tensor (2g - 1)}}{\frac{1}{2}} .
\end{equation*}
Since $h^2 = \gamma^2 = \frac{1}{2} p$ we get
\begin{equation*}
\GW{p, \gamma^{\tensor (2g + 1)}}{\frac{1}{2}} = \left( - \frac{1}{4} \right)^g \GW{p, \gamma}{\frac{1}{2}} 
\end{equation*}
by induction.  The invariant on the right side of this equality is easily seen to be $1$.  Indeed, $\overline{M}(\PP(1,1,2); 1, 1; \frac{1}{2})$ may be identified with $\PP(\Gamma(\PP(1,1,2), {\cal O}(1))) \cong \PP^1$.  The virtual dimension of $\overline{M}(\PP(1,1,2); 0, 1; \frac{1}{2})$ is also $1$, so we only need to solve the enumerative problem to compute $\GW{p, \gamma}{\frac{1}{2}}$.  If $(u, v) \in \PP^1$ is a point, then the condition that the corresponding curve interpolate the point $(x,y,z) \in \PP(1,1,2)$ is $ux + vy = 1$.  This has exactly one solution if $(x,y) \not= (0,0)$ so we conclude that $\GW{P, \gamma}{\frac{1}{2}} = 1$.

We have therefore proved
\begin{equation} \label{gw-calc}
\GW{p, \gamma^{\tensor (2g + 1)}}{\frac{1}{2}} = \left( - \frac{1}{4} \right)^g 
\end{equation}

\section{The virtual fundamental class}
\label{sect:interp}

\begin{lem-in-sect}
Let $C$ be a smooth orbifold curve.  Suppose there is a representable map $f : C \rightarrow \PP(1,1,2)$ of degree $\frac{1}{2}$.  Then $C$ has at most $1$ orbifold point.
\end{lem-in-sect}
\begin{proof}
In this case, $f^\ast {\cal O}(1)$ is a line bundle of degree $\frac{1}{2}$ on $C$.  The only such line bundles on $C$ are the ${\cal O}(P)$ where $P$ is an orbifold point of $C$.  Suppose $f^\ast {\cal O}(1) = {\cal O}(P)$ for a particular orbifold point $P$ and that $C$ has another orbifold point $Q \not= P$.  Let $\sigma$ be any nonzero section of ${\cal O}(1)$ over $\PP(1,1,2)$.  Then $f^\ast \sigma$ is a section of ${\cal O}_C(P)$, hence vanishes only $P$.  But this means $f(Q) \not= (0,0,1)$, which contradicts the representability of $f$.
\end{proof}

\begin{prop-in-sect} \label{modspisom}
There are isomorphisms
\begin{equation*}
\overline{M}(\PP(1,1,2); 0, 2g + 1; \frac{1}{2}) \cong \overline{M}(\PP(1,1,2); 0, 1; \frac{1}{2}) \times \overline{M}(B (\ZZ / 2 \ZZ); 2g + 2) 
\end{equation*}
\end{prop-in-sect}
\begin{proof}
  If $(f, C) \in \overline{M}(\PP(1,1,2); 0, 2g + 1; \frac{1}{2})$, then $C$ has a unique irreducible component $C_0$ with $\deg f \rest{C_0} = \frac{1}{2}$; all other components have degree $0$.  By the lemma, $C_0$ has exactly $1$ orbifold point.  The remaining orbifold points must lie on a component that is attached at the unique orbifold point of $C$.  Thus every point of $\overline{M}(\PP(1,1,2); 0, 2g + 1; \frac{1}{2})$ lies in the image of the gluing map
\begin{equation*}
\iota : \overline{M}(\PP(1,1,2); 0, 1; \frac{1}{2}) \times \overline{M}(B (\ZZ / 2 \ZZ); 2g + 2) \rightarrow \overline{M}(\PP(1,1,2); 0, 2g + 1) 
\end{equation*}
that attaches the marked point from the first component to the first marked point from the second component.

This is a closed embedding, so to complete the proof, we must show that the image of this map is in open in $\overline{M}(B (\ZZ / 2 \ZZ); 2g + 2)$.  Consider a first-order deformation $(C', f')$ of $(C, f)$.  Let $C_1$ be the contracted component of $C$ and let $C_0$ be the component of positive degree.  If $(C', f')$ were not in the  image of $\iota$, then $C'$ would be a first-order smoothing of $C$.  But then consider the map $N_{C_1 / C'} \rightarrow (f \rest{C_1})^\ast T_{(0,0,1)} \PP(1,1,2)$.  If $P$ is the point of attachment between $C_0$ and $C_1$, then $N_{C_1 / C'} \rest{P}$ is spanned by $T_P C_0$.  Moreover, $C_0$ meets $(0,0,1)$ transversally (since $f \rest{C_0}$ has degree $\frac{1}{2}$), which implies that the map $N_{C_1 / C'} \rightarrow (f \rest{C_1})^\ast T_{(0,0,1)} \PP(1,1,2)$ is nonzero at $P$.

On the other hand, $N_{C_1 / C'} \cong {\cal O}_{C_1}(-P)$, and $(f \rest{C_1})^\ast T_{(0,0,1)} \PP(1,1,2) \cong (f \rest{C_1})^\ast (\rho_1 \oplus \rho_1)$ because $f$ contracts $C_1$ onto the point $(0,0,1)$ and $T_{(0,0,1)} \PP(1,1,2) \cong \rho_1 \oplus \rho_1$.  Thus we obtain a pair of sections of $\rho_1 \tensor {\cal O}_{C_1}(P)$, at least one of which does not vanish at $P$.

Let $\pi : C_1 \rightarrow \overline{C}_1$ be the coarse moduli space.  Then we get a section of $\pi_\ast(\rho_1 \tensor {\cal O}_{C_1}(P))$ that is not everywhere zero.  But $\pi_\ast(\rho_1 \tensor {\cal O}_{C_1}(P)) = {\cal O}_{\overline{C}_1}(-g)$ where $2g + 2$ is the number of orbifold points on $C_1$.  By stability of $(C, f)$, $g > 0$, so all sections of $\pi_\ast(\rho_1 \tensor {\cal O}_{C_1}(P))$ vanish.  This contradicts the nonvanishing of the section at $P$.
\end{proof}

Now that we know how the moduli space looks, we must determine the virtual fundamental class.  We use the deformation--obstruction sequence,
\begin{equation*}
\Def(C) \rightarrow \Obs(f) \rightarrow \Obs(C,f) \rightarrow \Obs(C) = 0 .
\end{equation*}
We know that $\Obs(C,f)$ is a vector bundle because $\overline{M}(\PP(1,1,2); 0, 2g + 1; \frac{1}{2})$ is smooth.   The virtual fundamental class is the top Chern class of this vector bundle.

$\Obs(f)$ is the relative obstruction space for the map 
\begin{equation*}
\overline{M}(\PP(1,1,2); 0, 2g + 1; \frac{1}{2}) \rightarrow {\frk M}(B(\ZZ / 2 \ZZ); 2g + 1).
\end{equation*}
If $(C, f)$ is a curve in $\overline{M}(\PP(1,1,2); 0, 2g + 1; \frac{1}{2})$ then we have just seen that $C$ is the union of two curves, $C_0$ and $C_1$, along an orbifold point, with $\deg(f \rest{C_0}) = \frac{1}{2}$ and $\deg(f \rest{C_1}) = 0$.  It is clear that any deformation of $C$ that is trivial near the node will extend to a deformation of $(C, f)$ --- indeed, $C_0$ is rigid and $C_1$ is contracted by $f$.  Thus, the image of $\Def(C) \rightarrow \Obs(f)$ is the space of deformations of the node.  If we name the nodal point $P$, then the deformations of the node are parameterized by $\pi_\ast (T_P C_0 \tensor T_P C_1)$, so we have an exact sequence on $\overline{M}(\PP(1,1,2); 0, 2g + 1; \frac{1}{2})$,
\begin{equation*}
0 \rightarrow \pi_\ast (T_P C_0 \tensor T_P C_1) \rightarrow \Obs(f) \rightarrow \Obs(C,f) \rightarrow 0 .
\end{equation*}

Explicitly, $\Obs(f) = R^1 \pi_\ast f^\ast T \PP(1,1,2)$, where $f : C \rightarrow \PP(1,1,2)$ is the universal map.  Tensoring the normalization sequence for the node $P$ with $f^\ast T \PP(1,1,2)$ and taking cohomology, we obtain
\begin{equation*}
H^0(T \rest{P}) \rightarrow H^1(T) \rightarrow H^1(T \rest{C_0}) \oplus H^1(T \rest{C_1}) \rightarrow H^1(T \rest{P}) = 0 ,
\end{equation*}
writing $T = f^\ast T \PP(1,1,2)$.  Note that $H^0(T \rest{P}) = 0$ since $P$ is an orbifold point and $T \rest{P} \cong \rho_1 \oplus \rho_1$ has no invariant sections.

We can also calculate $H^1(T \rest{C_0}) = 0$ using the Euler sequence, which pulls back to
\begin{equation*}
0 \rightarrow {\cal O} \rightarrow {\cal O}(P) \oplus {\cal O}(P) \oplus {\cal O}(2P) \rightarrow T \rest{C_0} \rightarrow 0
\end{equation*}
since $f \rest{C_0}$ has degree $\frac{1}{2}$ and $f^\ast {\cal O}(1) = {\cal O}(P)$.  Pushing this sequence forward to the coarse moduli space via $q : C_0 \rightarrow \overline{C}_0$ (note $q_\ast$ is exact) gives
\begin{equation*}
0 \rightarrow {\cal O}_{\overline{C}_0} \rightarrow {\cal O}_{\overline{C}_0} \oplus {\cal O}_{\overline{C}_0} \oplus {\cal O}_{\overline{C}_0}(q(P)) \rightarrow \pi_\ast T \rightarrow 0 .
\end{equation*}
Now taking cohomology and noting that $H^1({\cal O}_{\overline{C}_0}) = H^1({\cal O}_{\overline{C}_0}(\pi(P))) = H^2({\cal O}_{\overline{C}_0}) = 0$, we deduce that $H^1(T \rest{C_0}) = 0$ from the long exact sequence.

It now follows that $\Obs(f) = H^1(T \rest{C_1})$.  But, as already remarked, $f \rest{C_1}$ factors through the orbifold point of $\PP(1,1,2)$, so $T \rest{C_1}$ is the pullback of the tangent bundle at this point, which is $\rho_1 \oplus \rho_1$.  Thus,
\begin{equation*}
\Obs(f) = R^1 \pi_\ast (\rho_1 \oplus \rho_1) \cong {\bf E}^\vee \oplus {\bf E}^\vee
\end{equation*}
where ${\bf E}$ is the Hodge bundle defined in the introduction.

We therefore have an exact sequence,
\begin{equation*}
0 \rightarrow \pi_\ast(T_P C_0 \tensor T_P C_1) \rightarrow {\bf E}^\vee \oplus {\bf E}^\vee \rightarrow \Obs(C,f) \rightarrow 0 .
\end{equation*}

Now, consider the cartesian diagram
\begin{equation*} \xymatrix{
e^{-1}(p) \ar[r] \ar[d] & \overline{M}(\PP(1,1,2); 1, 2g + 1; \frac{1}{2}) \ar[d]^e \\
p \ar[r]^i & \PP(1,1,2) .
}\end{equation*}

Under the identification of Proposition~\ref{modspisom}, $e$ factors through the evaluation map on $\overline{M}(\PP(1,1,2); 1, 1; \frac{1}{2})$.  Thus, $e^{-1}(p)$ may be identified with $\overline{M}(B(\ZZ / 2 \ZZ); 2g + 2)$.  We have $i^{-1}(T_P C_0) \cong \rho_1$, we get the exact sequence,
\begin{equation*}
0 \rightarrow L_1^\vee \rightarrow {\bf E}^\vee \oplus {\bf E}^\vee \rightarrow  i^\ast \Obs(C, f) \rightarrow 0 .
\end{equation*}
Now,
\begin{equation*}
\GW{p, \gamma^{\tensor (2g + 1)}}{\frac{1}{2}} = \int i^! [ \overline{M}(\PP(1,1,2); 1, 2g + 1; \frac{1}{2}) ]^{\rm vir} = \int_{\overline{M}(B(\ZZ / 2 \ZZ); 2g + 2)} \frac{c({\bf E}^\vee)^2}{c(L_1^\vee)} .
\end{equation*}
But we have also seen in Section~\ref{sect:eval} that
\begin{equation*}
\GW{p, \gamma^{\tensor (2g + 1)}}{\frac{1}{2}} = \left( - \frac{1}{4} \right)^g 
\end{equation*}
and this completes the proof of the theorem.

\bibliographystyle{plain}
\bibliography{thesis}

\begin{thebibliography}{1}

\bibitem{AGV}
Dan Abramovich, Tom Graber, and Angelo Vistoli.
\newblock {G}romov--{W}itten theory of {D}eligne--{M}umford stacks.
\newblock math.AG/0603151, 2006.

\bibitem{AV}
Dan Abramovich and Angelo Vistoli.
\newblock Compactifying the space of stable maps.
\newblock {\em J. Amer. Math. Soc.}, 15(1):27--75 (electronic), 2002.

\bibitem{JW}
Jonathan Wise.
\newblock The genus zero {G}romov--{W}itten invariants of
  $[\operatorname{Sym}^2 {P}^2]$, 2007.

\end{thebibliography}

\end{document}